\pdfoutput=1
\documentclass[a4paper,11pt]{article}
\usepackage[utf8]{inputenc}
\usepackage[T1]{fontenc}
\usepackage[english]{babel}
\usepackage{lmodern,microtype}
\usepackage{cite} % order citations
\usepackage{amsmath,amssymb,amsthm,mathtools,tikz,enumitem}
\setlist[enumerate]{label=\textit{(\roman*)},ref=\textit{(\roman*)}}
\usepackage[left=3.5cm,right=3.5cm,top=2cm,bottom=2cm]{geometry}
\usepackage{changepage} % for \begin{adjustwidth}
\usepackage[%hidelinks,
colorlinks,%
linkcolor=black!40!red,%
citecolor=green!40!black,%
urlcolor=red,%
pdfauthor={F\'abio Botler and Phablo F.S. Moura and T\'assio Naia},%
pdftitle={Seymour's Second Neighborhood Conjecture for orientations of (pseudo)random graphs},%
pdfkeywords={directed graph, random graph, neighborhood}%
]{hyperref}

\usepackage{tikz,mdframed} % figures + borders on theorem envs

\colorlet{tn/color/theorem}{green!60!black}
\colorlet{tn/color/claim}{gray!60!white}
\colorlet{tn/color/definition}{black}
\colorlet{tn/color/question}{orange!60!white}
\colorlet{tn/color/conjecture}{orange!60!white}
\colorlet{tn/color/defi}{red!50!black}

\theoremstyle{definition}

\newtheoremstyle{slthm}% <name>
{0pt}% <Space above> % was \bigskipamount
{0pt}% <Space below> % was 3pt
{\normalfont}% <Body font>
{}% <Indent amount>
{\bfseries\small}% <Theorem head font>
{.\!}% <Punctuation after theorem head>
{.6em}% <Space after theorem headi>
{\thmname{#1}\thmnumber{ #2}\thmnote{\,--\,\ignorespaces#3}}% <Theorem head spec (can be left empty, meaning `normal')>
\theoremstyle{slthm}

\newlength{\mythmbar}
\newlength{\myinsep}
\newlength{\myleftinsep}
\newlength{\mylen}
\newlength{\thmskip}
\mdfsetup{skipabove=\thmskip,skipbelow=\thmskip,usetwoside=false,nobreak}

\newcommand{\tndefmdstyle}[2]{%
  \mdfdefinestyle{tn/#1}{%
    linewidth=\mythmbar,%
    linecolor=tn/color/#2,%
    bottomline=false,%
    topline=false,%
    innerleftmargin=\myleftinsep,%
    leftmargin=-\myinsep,%
    innerrightmargin=\myinsep,%
    innertopmargin=1pt,%
    innerbottommargin=0pt,%
    rightmargin=\myinsep,%
    skipabove=\bigskipamount,%
    skipbelow=\bigskipamount,%
    userdefinedwidth=\mylen}}

\tndefmdstyle{base}{theorem}
\tndefmdstyle{claim}{claim}
\tndefmdstyle{definition}{definition}
\tndefmdstyle{question}{question}

\newmdtheoremenv[style=tn/base]{theorem}{Theorem}

\newcommand{\tnmdenv}[3]{%
  \newmdtheoremenv[style=#1]{#2}[theorem]{#3}}

\tnmdenv{tn/base}{lemma}{Lemma}
\tnmdenv{tn/base}{proposition}{Proposition}
\tnmdenv{tn/claim}{claim}{Claim}
\tnmdenv{tn/claim}{example}{Example}
\tnmdenv{tn/base}{corollary}{Corollary}
\tnmdenv{tn/base}{fact}{Fact}
\tnmdenv{tn/question}{question}{Question}
\tnmdenv{tn/question}{problem}{Problem}
\tnmdenv{tn/question}{conjecture}{Conjecture}
\tnmdenv{tn/definition}{definition}{Definition}
\tnmdenv{tn/definition}{remark}{Remark}
\tnmdenv{tn/definition}{hole}{Hole}

\newtheoremstyle{tnUnnumb}% <name>
{0pt}% <Space above> % was \bigskipamount
{0pt}% <Space below> % was 3pt
{\normalfont}% <Body font>
{}% <Indent amount>
{\bfseries\footnotesize}% <Theorem head font>
{.}% <Punctuation after theorem head>
{.5em}% <Space after theorem headi>
{\thmnote{#3}}% <Theorem head spec (can be left empty, meaning `normal')>

\theoremstyle{tnUnnumb}
\newmdtheoremenv[style=tn/base]{unnumtheorem}{Theorem}
\newmdtheoremenv[style=tn/base]{unnumconjecture}{Conjecture}
\newmdtheoremenv[style=tn/base]{unnumlemma}{Lemma}

\newcommand{\calb}{\mathcal{B}}

\newcommand{\cale}{\mathcal{E}}

\newcommand{\bbe}{\mathbb{E}}

\newcommand{\bbn}{\mathbb{N}}

\DeclareMathOperator{\bad}{BAD}

\DeclareMathOperator{\EE}{\bbe}
\DeclareMathOperator{\Var}{Var}

\DeclareSymbolFont{bbold}{U}{bbold}{m}{n}
\DeclareSymbolFontAlphabet{\mathbbold}{bbold}
\newcommand{\bbone}{\mathbbold{1}}

\newcommand{\good}{U}
\newcommand{\eps}{\varepsilon}
\newcommand{\ee}{\mathrm{e}}

\newcommand{\littleo}{\mathrm{o}}

\newcommand{\deq}{\coloneqq}
\renewcommand{\deq}{=}

\newcommand{\SNC}{\mathcal{S}}

\newcommand{\defi}[1]{%
  \emph{\color{red!60!black}#1}%
  }

\newcommand{\newproofenv}[3]{%
  \newenvironment{#1}[1][#2]{%
    \renewcommand{\qedsymbol}{\hfill{#3}}\begin{proof}[#2]%
  }{%
    \end{proof}\renewcommand{\qedsymbol}{\oldqed}%
  }%
}

\newproofenv{claimproof}{Proof}{$\diamond$}
\newproofenv{proofsketch}{Proof sketch}{$\bigtriangleup$}

\newenvironment{subproof}[1]{%
  \renewcommand{\qedsymbol}{\hfill{$\diamond$}}
  \begin{proof}[#1]%
  }{%
  \end{proof}\renewcommand{\qedsymbol}{\oldqed}%
}

\newcommand{\authormark}[1]{\textsuperscript{\,#1}}

\makeatletter
\renewcommand{\@biblabel}[1]{\textcolor{gray}{[\,}#1\textcolor{gray}{\,]}}
  \renewcommand\@openbib@code{% change the list parameters
    \setlength\labelwidth{2cm}%
    \setlength{\itemindent}{0cm}%
    \setlength{\leftmargin}{0cm}%
    \setlength\labelsep{.8em}%
    }
\makeatother

\begin{document} % --------------------------------------
% \linenumbers % turn line-numbering on LINENO-ON

\title{Seymour's Second Neighborhood Conjecture\\
  for orientations of (pseudo)random graphs}

\author{F\'abio Botler\authormark{1}%
  \and
  Phablo F.\,S. Moura\authormark{2}%
  \and T\'assio Naia\authormark{3}
}

\maketitle % --------------------------------------------

\newcommand{\showmark}[1]{%
  \hspace*{-1em}\makebox[1em][r]{\authormark{#1}\,}%
  \ignorespaces}

\vspace{-.6cm}
\begin{center}
\footnotesize
\showmark{1}
Programa de Engenharia de Sistemas e Computa\c c\~ao  \\
Instituto Alberto Luiz Coimbra
de P\'os-Gradua\c c\~ao e Pesquisa em Engenharia      \\
Universidade Federal do Rio de Janeiro, Brasil        \\
{\texttt{fbotler@cos.ufrj.br}}

\bigskip
\showmark{2}
Departamento de Ci\^encia da Computa\c c\~ao          \\
Instituto de Ci\^encias Exatas                        \\
Universidade Federal de Minas Gerais, Brasil          \\
{\texttt{phablo@dcc.ufmg.br}}

\bigskip
\showmark{3}
Departamento de Ci\^encia da Computa\c c\~ao          \\
Instituto de Matem\'atica e Estat\'\i stica           \\
Universidade de S\~ao Paulo, Brasil                   \\
{\texttt{tnaia@member.fsf.org}}
\end{center}

\begin{abstract}
  Seymour's Second Neighborhood Conjecture (SNC) states that
  every oriented graph contains a vertex
  whose second neighborhood is as large as its first neighborhood.
  We investigate the SNC for orientations of both binomial and pseudo random graphs,
  verifying the SNC asymptotically almost surely (a.a.s.)
  \begin{enumerate}
  \item for all orientations of $G(n,p)$ if $\limsup_{n\to\infty} p < 1/4$; and
  \item for a uniformly-random orientation of each weakly
    $(p,A\sqrt{np})$-bijumbled graph of order $n$ and
    density~$p$, where $p=\Omega(n^{-1/2})$ and $1-p = \Omega(n^{-1/6})$
    and $A>0$ is a universal constant independent of both $n$~and~$p$.
  \end{enumerate}
  We also show that a.a.s.\ the SNC holds
  for almost every orientation of~$G(n,p)$.
  More specifically, we prove that a.a.s.
  \begin{enumerate}[resume]
  \item
    for all $\eps > 0$ and $p=p(n)$
    with $\limsup_{n\to\infty} p \le 2/3-\eps$,
    every orientation of~$G(n,p)$
    with minimum outdegree~$\Omega_\eps(\sqrt{n})$
    satisfies the SNC; and
  \item
    for all $p=p(n)$, a random orientation of~$G(n,p)$ satisfies the SNC.
  \end{enumerate}
\end{abstract}

\section{Introduction}

An \defi{oriented graph} \(D\)
is a digraph obtained from a simple graph \(G\)
by assigning directions to its edges
(i.e., $D$ contains neither loops, nor parallel arcs,
nor directed cycles of length~\(2\));
we also call \(D\) an \defi{orientation} of \(G\).
Given \(i\in\bbn\),
the \defi{\(i\)-th neighborhood} of \(u\in V(D)\), denoted by \defi{\(N^i(u)\)},
is the set of vertices \(v\) for which a shortest directed path
from \(u\) to \(v\) has precisely \(i\)~arcs.
A~\defi{Seymour vertex} (see~\cite{2015:Seacrest})
is a vertex~\(u\) for which \(|N^2(u)|\geq |N^1(u)|\).
Seymour  conjectured
the following (see~\cite{1995:DeanLatka}).

\begin{conjecture}
  \label{conj:snc}
   Every oriented graph contains a Seymour vertex.
\end{conjecture}

Conjecture~\ref{conj:snc}, known as \emph{Seymour’s Second Neighborhood Conjecture} (SNC),
is a notorious open question
(see, e.g.,~\cite{2003:ChenShenYuster,2007:FidlerYuster,2012:Ghazal,2015:Seacrest}).
In particular, it was confirmed for tournaments
(orientations of~cliques)
by~Fisher~\cite{1996:Fisher}
and (with a purely combinatorial argument)
by~Havet and Thomass\'e~\cite{2000:HavetThomasee}; it was also studied
by Cohn, Godbole, Harkness and Zhang~\cite{2016:Cohn_etal}
for the random digraph model in which each ordered pair of vertices
is picked independently as an arc
with probability \(p<1/2\).
Throughout the paper,
we denote by \defi{$\SNC$} the set of graphs
$\{G:\text{all orientations of $G$ contain a Seymour vertex}\}$.

\smallskip

Our contribution comes from considering this combinatorial problem in a random
and pseudorandom setting
(see, e.g.,~\cite{ConlonGowers16,Schacht2016}).
More precisely, we
explore Conjecture~\ref{conj:snc} for orientations
of the binomial random graph~\defi{\(G(n,p)\)},
defined as the random graph with vertex set~$\{1,\ldots,n\}$ in which every
pair of vertices appears as an edge independently and with probability~$p$.

We say that an event $\cale$ holds \defi{asymptotically almost surely}
(a.a.s.)
if $\Pr[\cale]\to 1$ as $n\to\infty$.
If $G=G(n,p)$ is very sparse
(say, if $np \le (1-\eps) \ln n$ for large~$n$ and fixed~$\eps>0$), then
a.a.s.\ $G$ has an isolated vertex, which clearly is a Seymour vertex.
Our first result extends this observation
to much denser random graphs.

\begin{theorem}\label{t:snc-p<1/4}
  Let~$p\colon\bbn\to (0,1)$.
  If $\displaystyle\limsup_{n\to\infty} p < 1/4$,
  then a.a.s.\ $G(n,p)\in\SNC$.
\end{theorem}

If we impose restrictions on the orientations,
requiring, for example, somewhat large minimum outdegree,
the range of~$p$ can be further increased.

\begin{theorem}\label{t:gnp-min-outdeg}
  For every $\beta >0$, there exists $C=C(\beta)$
  such that the following holds for all~$p\colon\bbn\to (0,1)$.
  If $\displaystyle\limsup_{n\to\infty} p \le 2/3 -\beta$,
  then a.a.s.\ every orientation of $G(n,p)$
  with minimum degree at least $Cn^{1/2}$ contains a Seymour vertex.
\end{theorem}

For an even larger range of~$p$, we show that \emph{most} orientations
of~$G(n,p)$ contain a Seymour vertex;
i.e., Conjecture~\ref{conj:snc} holds
 for almost every (labeled) oriented graph.

\begin{theorem}\label{t:typical}
  Let~$p\colon\bbn\to (0,1)$
  and let~$G=G(n,p)$.
  If~\(D\) is chosen
  uniformly at random among the $2^{e(G)}$
  orientations of~$G$,
  then a.a.s.\ \(D\) has a Seymour vertex.
\end{theorem}

In fact, we prove a version of Theorem~\ref{t:typical}
in a more general setting,
namely
orientations of pseudorandom graphs
(see Section~\ref{s:pseudo-random}).

\begin{theorem}\label{t:typical-bij}
  There exists an absolute constant $C>1$ such that
  the following holds.
  Let~$G$ be a weakly~$(p,A\sqrt{np})$-bijumbled graph
  of order~$n$, where $\eps^3np^2 \ge A^2C$
  and $p < 1-15\sqrt{\eps}$.
  If~\(D\) is chosen
  uniformly at random among the $2^{e(G)}$
  possible orientations of~$G$,
  then a.a.s.\ \(D\) has a Seymour vertex.
\end{theorem}

This paper is organized as follows.
In Section~\ref{sec:wheel-free} we prove Conjecture~\ref{conj:snc}
for wheel-free graphs,
which implies the particular case of Theorem~\ref{t:snc-p<1/4}
when $n^2p^3\to 0$.
In Section~\ref{s:p-typical} we complete the proof
of Theorem~\ref{t:snc-p<1/4}
and prove Theorems~\ref{t:gnp-min-outdeg} and~\ref{t:typical}
using a set of standard properties
of \(G(n,p)\). These properties are collected in Definition~\ref{d:p-typical}
and Lemma~\ref{l:gnp-typical} (proved in Appendix~\ref{a:auxiliary}).
In Section~\ref{s:pseudo-random}, we introduce bijumbled graphs
and prove Theorem~\ref{t:typical-bij}.
We make a few further remarks in Section~\ref{sec:concluding-remarks}.

To avoid uninteresting technicalities, we omit floor and ceiling signs.
If $A$ and $B$ are sets of vertices, we denote by~\defi{$\vec e\,(A,B)$}
the number of arcs directed from $A$ to~$B$, by~\defi{$e(A,B)$} the number
of edges or arcs with one vertex in each set, and by~\defi{$e(A)$}
the number of edges or arcs with both vertices in~$A$.
The (underlying) \defi{neighborhood} of a vertex~$u$ is denoted by~\defi{$N(u)$},
and the \defi{codegree} of vertices $u,\,v$
is~$\defi{\ensuremath{\deg(u,v)}}=\bigl|N(u)\cap N(v)\bigr|$.

We remark that Theorem~\ref{t:snc-p<1/4} and a weaker version of Theorem~\ref{t:gnp-min-outdeg}
appeared in the extended abstracts~\cite{botler2022:seymour-dmd,botler2022:seymour-etc}.

\section{Wheel-free graphs}\label{sec:wheel-free}

A \defi{wheel} is  a graph obtained from a cycle $C$ by adding a
new vertex adjacent to all vertices in \(C\).
Firstly, we show that $G(n,p)$ is wheel-free when $p$ is small;
then prove that all wheel-free graphs satisfy Conjecture~\ref{conj:snc}.

\begin{lemma}\label{l:wheel}
  If $p\in(0,1)$ and $n^4p^6 < \eps / 16$,
  then $\Pr\bigl[\,\text{$G(n,p)$ is wheel-free}\,\bigr]\ge 1-\eps$.
\end{lemma}

\begin{proof}
  We can assume $\eps < 1$.
  Since $n^4p^6 < \eps/16$, we have that
  \begin{align}\label{e:np-wheel-bounds}
    np^2 < (\eps p^2/16)^{1/4} < 1/2.
  \end{align}
  Let $X=\sum_{k=4}^n X_k$, where $X_k$ denotes
  the number of wheels of order $k$ in $G(n,p)$.
  By the~linearity of~expectation,
  \begin{align}
    \EE X
    & = \sum_{k=4}^n \EE X_k
    = \sum_{k=4}^n \binom{n}{k}k\frac{(k-1)!}{2(k-1)}p^{2(k-1)} \nonumber\\
    & < n\sum_{k=4}^n (np^2)^{k-1}
    = n^4p^6\sum_{k=0}^{n-4} (np^2)^k
     \stackrel{\mathrm{G.S.}}{<}  \frac{n^4p^6}{1 - np^2}\label{e:GP}
    \stackrel{\eqref{e:np-wheel-bounds}}{<} 2n^4p^6 < \frac{\eps}{8} < \eps.
  \end{align}
  Where in~\eqref{e:GP} we use
  the formula $\sum_{i=0}^\infty r^i = (1-r)^{-1}$
  for the geometric series (G.S.) of ratio~$r = np^2<1$.
  Markov's inequality then yields
  \(
    \Pr[X \ge 1] \le \EE X < \eps.
  \)
\end{proof}

To show that every orientation of a wheel-free graph
has a Seymour vertex, we prove a slightly stronger result.
A digraph is \defi{locally cornering}
if the outneighborhood of each vertex induces
a digraph with a \defi{sink} (i.e., a vertex of outdegree~$0$).
The next proposition follows immediately by noting that,
in a locally cornering digraph, each
vertex of minimum outdegree is a
Seymour vertex.

\begin{proposition}\label{p:locally-cornering}
  Every locally cornering digraph has a Seymour vertex.
\end{proposition}

Lemma~\ref{l:wheel}~and Proposition~\ref{p:locally-cornering} immediately
yield the following corollary.

\begin{corollary}\label{cor:snc-small-p}
  If $p\in (0,1)$,
  and $n^4p^6 < \eps / 16$, then
  $\Pr\bigl[\,G(n,p) \in \SNC\,\bigr] \ge 1-\eps$.
\end{corollary}

\begin{proof}
  Note that every orientation of a wheel-free graph is locally cornering,
  since the (out)neighborhood of each vertex is a forest,
  and every oriented forest has a vertex with outdegree~0.
  Hence the result follows
  by Lemma~\ref{l:wheel} and~Proposition~\ref{p:locally-cornering}.
\end{proof}

\section{Typical graphs}\label{s:p-typical}

In this section we prove
that if  $\limsup_{n\to\infty} p< 1/4$, then
a.a.s.\ \(G(n,p)\in\SNC\).
We use a number of standard properties of~$G(n,p)$,
stated for convenience in Definition~\ref{d:p-typical}.

\begin{definition}\label{d:p-typical}
  Let $p\in (0,1)$.
  A graph $G$ of order~$n$ is \defi{$p$-typical} if
  the following hold.
  \begin{enumerate}

  \item \label{i:p-typical:1}
    For every~$X\subseteq V(G)$, we have
    \[
      \biggl|e(X) - \binom{|X|}{2}p\biggr|
      \le |X|\sqrt{3np(1-p) }
      + 2 n.
    \]

  \item \label{i:p-typical:sharp-XY}
      If $n'\ln n \le n''\le n$ or $n'=n''=n$,
      then all~$X,\,Y\subseteq V(G)$ with~$|X|,|Y|\le n'$
      satisfy
      \[
        \bigl|\,e(X,Y) -|X||Y|p\,\bigr|
        \le \sqrt{6n''p(1-p) |X||Y|}
        + 2n''.
    \]

\item \label{i:p-typical:3}
    For every $v\in V(G)$, we have
    \[
      |\deg(v) - np\,|
      \le \sqrt{6np(1-p)\ln n}
      + 2 \ln n.
    \]

  \item \label{i:p-typical:common-neigh}
    For every distinct $u,v\in V(G)$,
    we have
    \[
      \bigl|\,\deg(u,v) - (n-2)p^2\,\bigr|
      \le \sqrt{6np^2(1-p^2)\ln n}
      + 2\ln n.
    \]

  \end{enumerate}
\end{definition}

It can be shown, using standard Chernoff-type concentration inequalities,
that $G(n,p)$ is $p$-typical with high probability
(see Appendix~\ref{a:auxiliary}).

\begin{lemma}
  \label{l:gnp-typical}
  For every $p\colon\mathbb{N} \to (0,1)$,
  a.a.s.\ $G=G(n,p)$ is $p$-typical.
\end{lemma}

We also use the following property of graphs satisfying
Definition~\ref{d:p-typical}\,\ref{i:p-typical:1}.

\begin{lemma}
  \label{l:bad-degrees:2}
  Let $G$ be a graph of order~$n$ which satisfies
  Definition~\ref{d:p-typical}\,\ref{i:p-typical:1}, and fix~$a\in\bbn$.
  If $D$ is an orientation of~$G$
  and $B=\{v\in V(D):\deg_D^+(v)<a\}$, then
  \[
    |B|
    \le \frac{2}{p}(a-1) + 1 + \sqrt{\frac{12n(1-p)}{p}} + \frac{4n}{|B|p}.
  \]
\end{lemma}

\begin{proof}
  The lemma follows by multiplying all terms in the inequality below by~$2/|B|p$.
  \begin{equation*}
    |B|(a-1) \ge e(G[B])
    \stackrel{\text{\ref{d:p-typical}\ref{i:p-typical:1}}}{\geq}
    \binom{|B|}{2}p-
    |B|\sqrt{3np(1-p)}
    - 2n.\qedhere
  \end{equation*}
\end{proof}

\subsection{Proof of Theorem~\ref{t:snc-p<1/4}}

Let us outline the proof of Theorem~\ref{t:snc-p<1/4}.
Firstly, we find a vertex~$w$
whose outneighborhood contains many
vertices with large outdegree.
Then, we note that
$|N^1(w)|=O(np)$ and that $N^1(w)\cup N^2(w)$
cannot be too dense.
Finally, since many outneighbors of~$w$
have large outdegree,
we conclude that $N^1(w)\cup N^2(w)$ must contain at least
$2|N^1(w)|$~vertices, completing the proof.
This yields the following.

\begin{lemma}\label{l:snc-p-typical<1/4}
  Fix $0<\alpha<1/4$ and~$\eps > 0$.
  There is $n_1=n_1(\alpha,\eps)$
  such that $\SNC$ contains all  $p$-typical graphs of order $n$
  such that $n\ge n_1$ and~$\eps n^{-2/3}\le p \le 1/4 - \alpha$.
\end{lemma}

Lemma~\ref{l:snc-p-typical<1/4} is our last ingredient
for proving Theorem~\ref{t:snc-p<1/4}.
Indeed,
fix $\eps >0$, set~$\alpha=1/4-\limsup_{n\to\infty} p(n)$ and let
$n_0$ be~large enough so that $p(n)\le 1/4-\alpha$
and so that $G(n,p)$ is $p$-typical with probability
at least~$1-\eps$ for all $n\ge n_0$ (this is Lemma~\ref{l:gnp-typical}).
Now either $p < \eps n^{-2/3}$ or $\eps n^{-2/3}\le p(n)\le 1/4-\alpha$.
In the former case we use Corollary~\ref{cor:snc-small-p},
and in the latter case Lemma~\ref{l:snc-p-typical<1/4},
concluding either way that
\[\Pr\bigl[\,G(n,p)\in\SNC\,\bigr] \ge 1-\eps.\]

\begin{proof}[Proof of Lemma~\ref{l:snc-p-typical<1/4}]
  We may and shall assume (choosing $n_1$ accordingly) that
  $np$~is large enough whenever necessary.
  Fix an arbitrary orientation of~$G$.
  For simplicity, we write $G$ for both the oriented and
  underlying graphs. Let
  \[\displaystyle S \deq \{v\in V(G) : \deg^+(v)< (1-\alpha)np/2\}
  \]
  and $T=V(G)\setminus S$.
  Firstly, we show that $|T|\ge\alpha n/2$.
  This is clearly the case if $|S|< \alpha n$ (since $\alpha < 1/4 < 1-\alpha$);
  let us show that this also holds if $|S| \ge \alpha n$.
  Indeed, since~$ p \ge \eps n^{-2/3}$,
  from Lemma~\ref{l:bad-degrees:2}
  with \(a = (1-\alpha)np/2\) we obtain
  \begin{align*}
    |S|
    &  \le \frac{2(a -1)}{p}
      + 1 + \sqrt{\frac{12n(1-p)}{p}} + \frac{4n}{|S|p}
        < (1-\alpha)n +\littleo(n)
      < \left(1-\frac{\alpha}{2}\right)n.
  \end{align*}
Therefore $|T| = n - |S|\ge\alpha n/2$
as desired.
Recall that \(np\) is large and \(p\leq 1/4\).
Then \(\sqrt{3np(1-p)}  \ge 4/\alpha\),
and hence, from Definition~\ref{d:p-typical}\,\ref{i:p-typical:1} ,
we get
\begin{align}
  e\bigl(T\bigr)
  & \ge  \binom{|T|}{2}p
    -    |T|\sqrt{3np(1-p)}
    -    2n
    >    \binom{|T|}{2}p
    -    2|T|\sqrt{np}
    >    \frac{|T|^2p}{3},\label{e:size-S}
\end{align}
and therefore, by averaging, there exists $w\in T$ satisfying
\begin{align}\label{eq:degree-w-lb}
  \deg_{T}^+(w)
  & \ge   \frac{e\bigl(T\bigr)}{|T|}
    \stackrel{\eqref{e:size-S}}{\geq}  \frac{\alpha n p}{6}.
\end{align}

We next show that $w$ is a Seymour vertex.
Let $X=N_G^1(w)$ and $Y=N_G^{2}(w)$, and suppose, for a contradiction,
that \(|Y| < |X|\).
From~\ref{d:p-typical}\,\ref{i:p-typical:3} and $p+\alpha \le1/4$, we have
\begin{align}
  |X|
  & \leq np +\sqrt{6np\ln n} +2\ln n
    < n\left(p +\frac{\alpha}{2}\right) < \frac{n}{4} \le \frac{n}{2}(1-2\alpha -2p).
    \label{e:size-X}
\end{align}
Moreover,
\begin{align}
  |X|
  =    \deg^+(w) \leq np +\sqrt{6np\ln n} +2\ln n
  &    < 2np.
    \label{e:size-X-2}
\end{align}
Recall that $w\in T$
and let~$N\deq X\cap T$ be the set of outneighbors of $w$ in~$T$.
By the definition of $N$ and~\eqref{eq:degree-w-lb} we have
\begin{align}
  |N|
  &   \ge  \frac{\alpha np}{6}.\label{e:|N|}
\end{align}

Note that $\vec e\,(N,X)$ counts arcs induced by~$N$ precisely once
(as~$N\subseteq X$),
and if the arc $u\to v$ is counted by $\vec e\,(N,X)$,
then $v$ is a common neighbor of $w$
and~$u\in N$.
Hence, by Definition~\ref{d:p-typical}\,\ref{i:p-typical:common-neigh},
we have that
\[\vec e\,(N,X) +e(N) \leq |N|\bigl(np^2 +\sqrt{6np^2\ln n} + 2\ln n\bigr).\]
Since vertices in \(T\) (and hence in \(N\))
have at least \((1-\alpha)np/2\) outneighbors,
we have
\begin{align}
\vec e\,(N,Y)
  & \geq   |N|\frac{(1-\alpha)np}{2} - \vec e\, (N,X) -e(N)\nonumber \\
  &\geq    |N|\frac{(1-\alpha)np}{2} - |N|\bigl(np^2+\sqrt{6np^2\ln n} +2\ln n\bigr)
  \label{eq:eNY:lb}.
\end{align}

The following estimate will be useful.
\begin{claim}\label{cl:small-terms}
It holds that \(2\ln n + \sqrt{6np^2\ln n} + \sqrt{6|Y|np/|N|} = \littleo(np)\).
\end{claim}

\begin{claimproof}
  We prove that each term in the sum above is \(\littleo(n)\)
  when divided by~\(p\).
  Clearly, \(\sqrt{6np^2\ln n}/p = \littleo(n)\).
  Recall that \(p\ge \eps n^{-2/3}\)
  and thus \((2\ln n) /p = o(n)\).
  Also,
   \begin{equation*}%\label{e:second-term}
     \sqrt{\frac{|Y|6n}{|N|p}}
     \stackrel{\eqref{e:|N|}}{\le}   \sqrt{\frac{|Y|36}{\alpha p^2}}
     <   \sqrt{\frac{|X|36}{\alpha p^2}}
     \stackrel{\eqref{e:size-X-2}}{<}   \sqrt{\frac{72n}{\alpha p}}
     =     \littleo(n).\qedhere
   \end{equation*}
\end{claimproof}

We divide the remainder of the proof in two cases. Fix \(\gamma\in(1/2, 2/3)\).

\smallskip\noindent\textbf{Case 1.}
Suppose firstly that~\(p > n^{\gamma - 1}/2\).
Using Definition~\ref{d:p-typical}\,\ref{i:p-typical:sharp-XY} we obtain
\begin{align}
  \label{eq:eNY:ub}
  \vec e\,(N,Y)
  \le |N||Y|p + \sqrt{6np |N||Y|} + 2n.
\end{align}
Thus, combining~\eqref{eq:eNY:lb} and \eqref{eq:eNY:ub}, we have
\begin{align}
\label{e:Y}
  \frac{(1-\alpha)np}{2} - (np^2+\sqrt{6np^2\ln n} + 2\ln n)
   & \le   |Y|p + \sqrt{\frac{6np|Y|}{|N|}} +\frac{2n}{|N|}.
\end{align}
Also note that since \(p > n^{\gamma - 1}/2\) and \(\gamma > 1/2\),
  we can estimate
\begin{align}
\label{e:bad-term-p-big}
  \frac{2n}{|N|p}
  \stackrel{\eqref{e:|N|}}{\leq}  \frac{12}{\alpha p^2}
  < \frac{24}{\alpha n^{2\gamma -2}}
  & = \littleo(n).
\end{align}
Finally, we conclude that $w$ is a Seymour vertex, since~\eqref{e:Y} becomes
\begin{align*}
 |Y|
  & \ge
    \frac{(1-\alpha -2p)n}{2} - \sqrt{6n\ln n} -\sqrt{\frac{6n|Y|}{|N|p}} -\frac{2n}{|N|p} -\frac{2\ln n}{p} \\
  & \stackrel{(\star)}{\ge}
    \frac{n}{2}\left(1 -2\alpha -2p\right)
   \stackrel{\eqref{e:size-X}}{>} |X|,
\end{align*}
where inequality~$(\star)$ follows from Claim~\ref{cl:small-terms} and~\eqref{e:bad-term-p-big}.

\smallskip\noindent\textbf{Case 2.}
Suppose now that \(p\leq n^{\gamma -1}/2\).
In this case \eqref{e:size-X-2}
implies~\(|X| \leq n^{\gamma}\).
Since \(N\subseteq X\) and \(|Y| < |X|\),
Definition~\ref{d:p-typical}\,\ref{i:p-typical:sharp-XY} (with $n'=n^\gamma$ and $n''=n^\gamma\ln n$) yields
\begin{align}
\vec e\,(N,Y)
   & \le |N||Y|p + \sqrt{6 (n^\gamma\ln n) p|N||Y|} + 2n^\gamma\ln n \nonumber\\
   &  <  |N||Y|p + \sqrt{6np |N||Y|} + 2n^\gamma\ln n.
       \label{eq:eNY:ub:p-small}
\end{align}

Now, from~\eqref{eq:eNY:lb} and~\eqref{eq:eNY:ub:p-small}, we obtain the following inequality,
which is analogous to~\eqref{e:Y}, but with the term \(2n/|N|\) replaced by \(2n^\gamma\ln n/|N|\).
\begin{equation}\label{e:Y:p-small}
  \frac{(1-\alpha)np}{2} - \bigl(p^2n+\sqrt{6np^2\ln n} + 2\ln n\bigr)
   \le   |Y|p + \sqrt{\frac{6np|Y|}{|N|}} +\frac{2n^\gamma\ln n}{|N|}.
\end{equation}
We claim that \(2n^\gamma\ln n/|N| = \littleo(np)\).
Indeed, since \(p\ge \eps n^{-2/3}\) and $\gamma < 2/3$, we have
\begin{equation}\label{e:last-term-estimate}
  \frac{2n^\gamma\ln n}{|N|p}
  \stackrel{\eqref{e:|N|}}{\leq}
  \frac{12n^\gamma\ln n}{\alpha np^2}
  =      \frac{12n^{\gamma-1}\ln n}{\alpha p^2}
  \le    \frac{12n^{\gamma+1/3}\ln n}{\alpha \eps^2}
  =      \littleo(n),
\end{equation}
We complete the proof of Case~2 by solving~\eqref{e:Y:p-small}
for~\(|Y|\) as in~Case~1 (using Claim~\ref{cl:small-terms}
and~\eqref{e:last-term-estimate} to estimate~\(2n^\gamma\ln n/|N|\)\,).
\end{proof}

  \subsection{Proof of Theorem~\ref{t:typical}}

  We are now in a position to prove Theorem~\ref{t:typical},
  which we restate for convenience.

\begin{unnumtheorem}[Theorem~\ref{t:typical}]
  Let~$p\colon\bbn\to (0,1)$,
  and let~$G=G(n,p)$.
  If~\(D\) is chosen
  uniformly at random among the $2^{e(G)}$
  orientations of~$G$,
  then a.a.s.\ \(D\) has a Seymour vertex.
\end{unnumtheorem}

\begin{proof}[Proof of Theorem~\ref{t:typical}]
  Let $G=G(n,p)$.
  If $p<1/5$, then $\Pr[\,G\in \SNC\,]= 1-\littleo(1)$ by Theorem~\ref{t:snc-p<1/4}.
  On the other hand, if $p\ge 1/5$,
  then standard concentration results for binomial
  random variables (e.g., Chernoff-type bounds)
  yield that every ordered pair $(u,v)$ of distinct vertices
  of $G$ satisfies, say $\deg(u,v)\ge n/50$, and hence with probability $1-\littleo(1)$
  every such pair is joined by a directed path of length~$2$.
  This is because building a random orientation of $G(n,p)$ is equivalent to first choosing
  which edges are present and then choosing the orientation of each edge uniformly at random,
  with choices mutually independent for each edge.
  In other words, with probability $1-\littleo(1)$,
  for all $u\in V(G)$ we have $V(G)=\{u\}\cup N^1(u)\cup N^2(u)$.
  Finally, by averaging outdegrees,
  we can find a vertex~$z\in V(D)$ with outdegree at most~$(n-1)/2$,
  because~$\sum_{v\in V(D)}\deg^+(v)=e(G) \le n(n-1)/2$.
  Such \(z\) is a Seymour vertex as desired.
\end{proof}

\subsection{Orientations with large minimum outdegree}
\label{s:minimum-degree}

Our last result in this section yields yet another class of
orientations of $p$-typical graphs which must always contain
a Seymour vertex.
In fact, we consider a larger class of underlying graphs,
showing that
if a graph~$G$ satisfies
items
\ref{i:p-typical:1}~and~\ref{i:p-typical:sharp-XY}
of Definition~\ref{d:p-typical}, then
every orientation~$D$ of $G$
with minimum outdegree~$\delta^+(D)=\Omega(n^{1/2})$
contains a Seymour vertex.
This may be useful towards extending the range of~$p$
for which a.a.s.\ $G(n,p)\in\SNC$.

\begin{lemma}\label{l:min-degree}
  Fix $\beta > 0$.
  There exist a constant~$C=C(\beta)$
  and $n_0=n_0(\beta)$
  such that the following holds for all $n\ge n_0$ and $p\le 2/3 -\beta$.
  If $G$ is a graph of order $n$
  that satisfies items~\ref{i:p-typical:1}
  and~\ref{i:p-typical:sharp-XY}
  of Definition~\ref{d:p-typical},
  then every orientation $D$ of~$G$
  for which~$\delta^+(D)\ge C n^{1/2}$ has a Seymour vertex.
\end{lemma}

Note that
Lemma~\ref{l:min-degree} and Lemma~\ref{l:gnp-typical}
immediately imply Theorem~\ref{t:gnp-min-outdeg}.

\begin{proof}[Proof of Lemma~\ref{l:min-degree}]
  Since $(1-3p/2)\ge 3\beta/2$, we may fix $C \ge 4$ so that
  \[
    \Bigl(1-\frac{3p}{2}\Bigr)C - \Bigl(\sqrt{3p(1-p)} + \sqrt{6p(1-p)}\Bigr)
    \geq \frac{3\beta C}{2} - 4 \ge 1.
  \]
  Fix \(v\in V(D)\) with $\deg^+(v)=\delta^+(D)$,
  let \(X \deq N^1(v)\) and \(Y \deq N^2(v)\).
  We shall prove that \(|X| \leq |Y|\).
  Suppose to the contrary that \(|Y| < |X|\).
  By Definition~\ref{d:p-typical}\,\ref{i:p-typical:1},
  \begin{align}\label{e:e(X,Y)-lower}
    \vec e\,(X,Y)
    =  \sum_{a\in X}\deg^+(a) - e(X)
    & \ge |X|^2 - \biggl(\frac{|X|^2p}{2} + |X|\sqrt{3np(1-p)} + 2n\biggr) \nonumber\\
    & = \left(1-\frac{p}{2}\right)|X|^2 - \Bigl(|X|\sqrt{3np(1-p)} + 2n\Bigr),
  \end{align}
  and by Definition~\ref{d:p-typical}\,\ref{i:p-typical:sharp-XY} (with $n'=n''=n$)
  we have
  \begin{align}\label{e:e(X,Y)-upper}
    \vec e\,(X,Y)
    \le e(X,Y)
    & \le |X||Y|p + \sqrt{6np(1-p)|X||Y|} + 2n \nonumber\\
    & < |X|^2p + |X|\sqrt{6np(1-p)} + 2n.
  \end{align}
  Since $|X|\ge Cn^{1/2} \ge n^{1/2}$,
  combining \eqref{e:e(X,Y)-lower}~and~\eqref{e:e(X,Y)-upper}
  yields the following contradiction.
  \begin{align*}
  4n   & > \Bigl(1-\frac{3p}{2}\Bigr)|X|^2 - |X|\Bigl(\sqrt{3np(1-p)} + \sqrt{6np(1-p)}\Bigr) \\
   & \ge Cn\left( \Bigl(1-\frac{3p}{2}\Bigr)C - \Bigl(\sqrt{3p(1-p)} + \sqrt{6p(1-p)}\Bigr)\right) \geq 4n.
          \qedhere
  \end{align*}
\end{proof}

\section{Typical orientations of bijumbled graphs}
\label{s:pseudo-random}

In this section, we focus on a well-known class of pseudorandom
graphs (that is, deterministic graphs which embody many
properties of $G(n,p)$\,), and argue that almost all of their orientations
 contain a Seymour vertex.
The following results concern graphs of order $n$ and density~$p$,
where $C n^{-1/2}\le p\le 1- \eps$, and $C = C(\eps) >0$ depends only on the constant~$\eps >0$.

\begin{definition}[$(p,\alpha)$-bijumbled]
  \label{def:(p,alpha)-bijumbled}
  Let~$p$ and~$\alpha$ be given.  We say that a graph~$G$ of order~$n$
  is \defi{weakly $(p,\alpha)$-bijumbled} if, for all~$U$,
  $W\subset V(G)$ with $U\cap W=\emptyset$
  and~$1\leq|U|\leq|W|\leq np|U|$, we have
  \begin{equation}
    \label{eq:(p,alpha)-bijumbled_def}
    \big|e(U,W)-p|U||W|\big|\leq\alpha\sqrt{|U||W|}.
  \end{equation}
  If~\eqref{eq:(p,alpha)-bijumbled_def}
  holds for all disjoint~$U$, $W\subset V(G)$, then we say that~$G$ is
  \defi{$(p,\alpha)$-bijumbled}.
\end{definition}

We note that the random graph is a.a.s.\ bijumbled.

\begin{theorem}[Lemma 3.8 in \cite{haxell95ramsey}]
  \label{t:gnp-is-bijumled}
  For any $p:\bbn\to (0,1]$,
  the random graph $G(n,p)$ is a.a.s.\
  weakly $(p, A\sqrt{np})$-bijumbled for a certain absolute
  constant~$A\le \ee^2\sqrt{6}$.
\end{theorem}

In what follows, $A$ shall always denote the constant from Theorem~\ref{t:gnp-is-bijumled}.
A simple double-counting argument shows the following.

\begin{fact}
  \label{fact:jumbled}
  If~$G$ is weakly $(p,\alpha)$-bijumbled, then for
  every~$U\subset V(G)$ we have
  \begin{equation}
    \label{eq:(p,alpha)-single}
    \left|
      e\bigl(G[U]\bigr) - p \binom{|U|}{2}
    \right| \leq \alpha|U|.
  \end{equation}
\end{fact}

We also use the following result,
whose simple proof we include for completeness.

\begin{lemma}\label{l:wbij-props}
  There exists a universal constant~$C>1$ such that
  if $A\ge 2$ and $\eps,\,p\in(0,1)$  are such that $\eps^3np^2 \ge A^2C$,
  then every weakly~$(p,A\sqrt{np})$-bijumbled graph~$G$ of order~$n$ satisfies
  the following properties.
  \begin{enumerate}
  \item\label{l:bij-props:deg}
    $\bigl|\,\{v\in V(G) :
    |\deg(v)-np| > \eps np\}\,\bigr| \le \eps n$.
  \item\label{l:bij-props:codeg}
    $\bigl|\,\{(u,v)\in V(G)^2 :
    \deg(u,v) \le (1 - \eps)np^2\}\,\bigr| \le \eps n^2$.
  \item\label{l:bij-props:orient}
    For every orientation of $G$ and every integer $d$,
    we have
    \[
      \bigl|\,\{v\in V(G) : \deg^+(v) < d\}\,\bigr|
      \leq 2\frac{d-1}{p} + 2A\sqrt{\frac{n}{p}} +1
    \]
  \end{enumerate}
\end{lemma}

  \begin{proof}
    Let \(G\) be as in the statement.
    We may and shall assume that $C$ is large enough so that the required inequalities
    hold.
    Throughout this proof, $W$
    denotes the set of vertices with degree strictly below~$(1-2\eps/3)np$.
    Firstly, we prove~\ref{l:bij-props:deg}.
    We claim that $|W|<\eps n/2$.
    Indeed, suppose the contrary and
    consider a subset \(W'\subseteq W\) of size precisely \(\eps n/2\).
    By Fact~\ref{fact:jumbled}, we have
    \begin{align}
      e(W')
      &  \ge   p\frac{(\eps n/2)^2}{3} - A\sqrt{np}(\eps n/2)
         =      p\frac{(\eps n/2)^2}{3}\left(1 - \sqrt{\frac{36A^2}{\eps^2np}}\right)
      \nonumber\\
      & >      p\frac{(\eps n/2)^2}{4}
        =      \frac{\sqrt{2}}{16}
        \,     \sqrt{\frac{\eps^3np}{1-\eps/2}}
        \,     \sqrt{\vphantom{\frac{\eps^3np}{1-\eps/2}}\frac{\eps}{2} (1-\eps/2)n^3p}
      \nonumber\\
      & \ge    A\sqrt{np\frac{\eps n}{2}(1-\eps/2)n}
        =      A\sqrt{np|W'|(n-|W'|)}.
        \label{e:e(W')}
    \shortintertext{Now, note that $|V(G)\setminus W'| < n  < A^2Cn/(\eps^2p) \le \eps n^2p = np|W'|$, but}
    e\bigl(W', V(G)\setminus W'\bigr)
    & <   |W'| \cdot (1-2\eps/3)np -2e(W')\nonumber\\
    & <   |W'| \cdot (1-\eps/2)np -2e(W') \nonumber\\
    & =   p|W'|(n-|W'|) -2e(W')           \nonumber\\
    & \stackrel{\eqref{e:e(W')}}{\le} p|W'|(n-|W'|) - A\sqrt{np|W'|(n-|W'|)}, \nonumber
    \end{align}
    which contradicts the weak bijumbledness of \(G\).

    Similarly, we show that the set $Z$ of vertices
    having degree strictly greater than $(1+2\eps/3)pn$
    satisfies $|Z| < \eps n/2$, which together with the argument
    above proves~\ref{l:bij-props:deg}.
    More precisely,
    suppose $|Z|\ge \eps n / 2$,
    fix $Z'\subseteq Z$ with $|Z'|=\eps n/2$.
    We claim that
    $A\sqrt{np}|Z'|$ and~$A\sqrt{np|Z'|(n-|Z'|)}$
    are both small (constant) fractions of~$p|Z'|^2$.
    Indeed, as $|Z'|^2  < |Z'|(n-|Z'|) < |Z'|n $, it follows that
    \begin{align*}
      \frac{A\sqrt{np}|Z'|}{p|Z'|^2}
      & < \frac{A\sqrt{np|Z'|(n-|Z'|)}}{p|Z'|^2}
        < \frac{A\sqrt{n^2p|Z'|}}{p|Z'|^2}\\
      & =   \sqrt{\frac{A^2n^2}{p|Z'|^3}}
        =   \sqrt{\frac{A^2n^2}{p(\eps n/2)^3}}
        \stackrel{(\sharp)}{\le} \sqrt{\frac{8p}{C}},
    \end{align*}
    where \((\sharp)\) is due to~$\eps^3np^2\ge CA^2$.
    Fact~\ref{fact:jumbled} and the previous inequalities imply
    \begin{align*}
      e(Z')
      & <   \frac{p|Z'|^2}{2} + A\sqrt{np}|Z'|
      \\
      & <  p|Z'|^2\left(\frac{1}{2}+\sqrt{\frac{8p}{C}}\right) \\
      & <  p|Z'|^2\left(\frac{1}{2}+\sqrt{\frac{32p}{C}}\right) -A\sqrt{np|Z'|(n-|Z'|)}
      \\
      & <  p|Z'|^2 -A\sqrt{np|Z'|(n-|Z'|)}.
    \end{align*}
    Analogously, we have $|V(G)\setminus Z'| < np|Z'|$, but
    \begin{align*}
      e(Z',V(G)\setminus Z')
      & \ge   (1+2\eps/3)np|Z'| - 2e(Z') \\
      & \ge   p|Z'|\bigl(n-|Z'|\bigr)
        + \left(\frac{1}{2}+\frac{2}{3}\right)\eps np|Z'| - 2e(Z') \\
      & >   p|Z'|\bigl(n-|Z'|\bigr) + 2p|Z'|^2 - 2e(Z') \\
      & >   p|Z'|\bigl(n-|Z'|\bigr) + A\sqrt{np|Z'|(n-|Z'|)},
    \end{align*}
    which is again a contradiction to
    Definition~\ref{def:(p,alpha)-bijumbled}.
    This concludes the proof of~\ref{l:bij-props:deg}.

    We next prove~\ref{l:bij-props:codeg}.
    For each $u\in V(G)$, let $B(u)$ be the set of vertices
    that have fewer than $(1-\eps)np^2$ common neighbors with~$u$.
    By definition,
    for any vertex \(u\) and set \(B'\subseteq B(u)\) we have
    \(e\bigl(N(u), B'\bigr) < (1-\eps)np^2\bigl|B'\bigr|\).
    We shall prove that $\bigl|B(u)\bigr|< \eps n/2$ for all~$u\in V(G)\setminus W$.
    Indeed,  suppose
    for a contradiction, that~$u\in V(G)\setminus W$ and~$|B(u)|\ge \eps n/2$.
    Let \(N'\subset N(u)\) be a set of size precisely \((1-2\eps/3)np\),
    and let \(B'\subseteq B(u)\) be a set of size precisely \(\eps n/2\).
    Since $\eps^3np^2 \ge A^2C$, we have
    \begin{align}
      \frac{\eps n p^2|B'|}{3}
      = \frac{\eps^2 n^2 p^2}{6}
      > \frac{1}{6}\sqrt{\frac{\eps^4n^4p^4(1-2\eps/3)}{2}}
      & > A\sqrt{np |N'||B'|}.\label{e:aux-2}
    \end{align}
    We claim that $|B'|\le np|N'|$.
    Indeed, $|N'|\le np\le \eps n^2p/2= np |B'|$
    because~$\eps n/2 > 1$,
    and~$|B'|=\eps n/2 \le A^2Cn/(3\eps^3)\le n^2p^2/3 \le np|N'|$
    because~$\eps^3np^2 \ge A^2C$ and~$\eps <1$.
    Hence, since $G$ is weakly bijumbled,
    we reach the following contradiction
    \begin{align*}
      p|N'||B'| - A\sqrt{np |N'||B'|}
           \le   e\bigl(N', B'\bigr)
            & <      (1-\eps)np^2\bigl|B'\bigr| \\
      & =  \left(1-\frac{2\eps}{3}\right)np^2|B'| -\frac{\eps n p^2|B'|}{3}\\
            &  \stackrel{\eqref{e:aux-2}}{<}    p|N'||B'| - A\sqrt{np |N'||B'|}.
    \end{align*}
    Hence~$\bigl|B(u)\bigr|<\eps n/2$ for all~$u\in V(G)\setminus W$.
    Note that if~$|N(u)\cap N(v)| < np^2(1-\eps)$
    for distinct~$u,\,v\in V(G)$,
    then either $u\in W$ or $v\in B(u)$.
    We conclude that there are at most $|W|n + n(\eps n/2) < \eps n^2$ such pairs,
    as desired.

    To prove~\ref{l:bij-props:orient},
    fix an orientation~$D$ of $G$
    and put~\(X\deq\{v\in V(G):\deg_{D}^+(v)< d\}\).
    Fact~\ref{fact:jumbled} then yields the desired inequality:
  \[
    |X|(d-1) \ge e(G[B])
    \ge
    \binom{|X|}{2}p-
    A\sqrt{np}|X|.\qedhere
  \]
  \end{proof}

\subsection{Almost all orientations of bijumbled graphs}
%\label{s:typical}

In this section we show that
almost every orientation of a weakly bijumbled graph
contains a Seymour vertex.

\begin{unnumtheorem}[Theorem~\ref{t:typical-bij}]
  There exists an absolute constant $C>1$ such that
  the following holds.
  Let~$G$ be a weakly~$(p,A\sqrt{np})$-bijumbled graph
  of order~$n$, where $\eps^3np^2 \ge A^2C$
  and $p < 1-15\sqrt{\eps}$.
  If~\(D\) is chosen
  uniformly at random among the $2^{e(G)}$
  possible orientations of~$G$,
  then a.a.s.\ \(D\) has a Seymour vertex.
\end{unnumtheorem}

\begin{proof}
  We may and shall assume that $A^2C$
  is larger than any given absolute constant.
    Let $V=V(G)$.
    For each $u\in V$,
    let $B(u)\deq\{v\in V:\deg(u,v)\le (1-\eps)np^2\}$.
    Also, let $\bad_1=\bigl\{u\in V:|B(u)|\ge \sqrt{\eps}n\bigr\}$.
    Lemma~\ref{l:wbij-props}\,\ref{l:bij-props:codeg} guarantees that $|\bad_1|\le \sqrt{\eps}n$
    and, by definition, $|B(u)|<\sqrt{\eps}n$
    for each $u\notin \bad_1$.

    Fix an arbitrary orientation of~$G$.
    For simplicity, we write~$G$ for both the oriented and underlying graphs.
    Let $\bad_2=\{v\in V(G): \deg^+(v) < 2\sqrt{\eps}np\}$.
    By~Lemma~\ref{l:wbij-props}\,\ref{l:bij-props:orient}, we must have
    \[
      |\bad_2|\le \frac{2(2\sqrt{\eps}np -1)}{p}+2A\sqrt{\frac{n}{p}} + 1 < 5\sqrt{\eps}n.
    \]
    Let $\bad=\bad_1\cup\bad_2$
    and put~$\good=V \setminus \bad$,
    and note that~$|\bad|\le 6\sqrt{\eps}n$.

   \begin{claim}\label{cl:deg(w)}
    There exists $w\in\good$
      such that
      \begin{equation*}\label{e:deg(w)}
        \deg_G^+(w) < n/2 - \sqrt{\eps}n.
      \end{equation*}
    \end{claim}
      \begin{subproof}{Proof}
              Recall that $p<1 -15\sqrt{\eps}$.
      Hence~$\eps < 15^{-2}<1$ and
    \begin{equation}\label{e:aux}
      \frac{(1+\eps)p}{2} + 6\sqrt{\eps}
      < \frac{(1+\eps)(1-15\sqrt{\eps})}{2} + 6\sqrt{\eps}
      < \frac{1-2\sqrt{\eps}}{2}
    \end{equation}
      Note also that $\eps^3np^2 \ge A^2C$ yields $A \le \sqrt{\eps^3np^2/C}$. Hence,
    \begin{equation}\label{e:anp-bound}
      A\sqrt{np}
      \le \eps np\sqrt{\frac{\eps p}{C}}
      < \frac{\eps np}{2}.
    \end{equation}
    By Fact~\ref{fact:jumbled}, we have
        \begin{align}\label{e:avg-deg}
          \frac{e(G[\good])}{|\good|}
          & \le  \frac{p}{|\good|}\binom{|\good|}{2}+A\sqrt{np}
          \le \frac{p|\good|}{2} + A\sqrt{np}
          \stackrel{\eqref{e:anp-bound}}{\le} (1+\eps)\frac{np}{2}.
    \end{align}
    Owing to~\eqref{e:avg-deg}, averaging the outdegrees of vertices in $\good$
    yields that some $w\in \good$
    satisfies $\deg_{G[\good]}^+(w)\le e\bigl(G[\good]\bigr)/|\good| < (1+\eps)np/2$.
    Hence,
    \begin{align}
      \deg_G^+(w)
      & \le  \deg_{G[\good]}^+(w) +|\bad|\nonumber\\
      & <    \frac{(1+\eps)np}{2} + 6\sqrt{\eps}n
       \stackrel{\eqref{e:aux}}{\le}    \frac{(1-2\sqrt{\eps})n}{2}.\qedhere
    \end{align}
  \end{subproof}

Note that since we picked an arbitrary orientation of \(G\),
the vertex \(w\) given by Claim~\ref{cl:deg(w)} exists for any such orientation.
  To conclude the proof, we next show that
  in a random orientation of~$G$
  almost surely every vertex in~$\good$
  is an~$(1-2\sqrt{\eps})$-king,
  where a vertex $v$ is said to be a \defi{$\lambda$-king}
  if the number of vertices $z$ for which
  there exists a directed path of length $2$ from~$v$ to $z$ is at least~$\lambda n$.

    \begin{claim}\label{cl:large-neigh}
      In a random orientation of $G$,
      a.a.s.\ for each $X\subseteq V(G)$
      with $|X|=2\sqrt{\eps}np$ we have $\bigl|N^1(X)\bigr|\ge (1-2\sqrt{\eps})n$,
      where $N^1(X)=\bigcup_{x\in X} N^1(x)$.
    \end{claim}
    \begin{claimproof}
      Note that for all $X,Y\subseteq V(G)$, there exist $X'\subseteq X$ and $Y'\subseteq Y$
      such that $X'\cap Y'=\emptyset$ and $|X'|= |X|/2$ and $|Y'|= |Y|/2$.
      Fix $X\subseteq V(G)$ with $|X|=2\sqrt{\eps}np$.
      If we choose $Y$ such that $|Y|=2\sqrt{\eps}n$,
      then~$|X'| \le |Y'| = \sqrt{\eps}n \leq \sqrt{\eps}n^2p^2 = np |X'|$
      because \(np^2 \geq A^2C/\eps^3 \geq 1\). Hence, as $G$ is weakly bijumbled,
      \begin{equation}\label{e:nedges}
        e(X,Y)\ge e(X',Y')\ge \frac{p|X||Y|}{4} -\frac{A\sqrt{np|X||Y|}}{2}.
      \end{equation}
      Let $\cale_X$ denote the `bad' event that $\bigl|N^1(X)\bigr|< (1-2\sqrt{\eps})n$,
      so $\cale_X$ occurs if and only if
      there exists $Y\subseteq V(G)$ with $|Y|=2\sqrt{\eps}n$ such that~$\vec e\, (X,Y)=0$.
      For any $X$ such that~$|X|=2\sqrt{\eps}np$, summing over all $Y$ of size~$2\sqrt{\eps}n$
      yields
      \begin{align}
        \Pr[\cale_X]
        & \le  \sum_{Y}2^{-e(X,Y)}
          \stackrel{\eqref{e:nedges}}{\le}
          \binom{n}{2\sqrt{\eps}n}
          \exp\Bigl( -(\ln 2)\bigl(\eps n^2p^2 -A\sqrt{\eps n^3p^2}\bigr)\Bigr)
         \nonumber\\
        & \le
          \exp\left(
          2n\sqrt{\eps}\ln\left(  \frac{\ee}{2\sqrt{\eps}}  \right)
          -(\ln 2)\eps n^2p^2\Bigl(1-\frac{\eps}{\sqrt{C}}\Bigr)
          \right)
          \nonumber\\
        &\le
          \exp\left(  2n\sqrt{\eps}\left(\frac{\ee}{2\sqrt{\eps}}  \right)
          -(\ln 2)\eps n^2p^2\Bigl( 1-\frac{\eps}{\sqrt{C}} \Bigr)\right)
          \nonumber\\
         &\le   \exp\bigl(-2(\ln 2)n\bigr)
      \label{e:Pr(cale_X)}
      \end{align}
      using that~$\eps np^2\ge A^2C\eps^{-2}\geq 12$ and that $\eps/\sqrt{C} \le C^{-1/2} < 1/2$
      because $\eps < 1$ and $C$ is a large constant.
      Taking a union bound over all~$X$ of size~$2\sqrt{\eps}np$,
      we see that no bad event occurs is with high probability,
      since
      \[
        \sum_X \Pr\bigl[ \cale_X \bigr]
        \stackrel{\eqref{e:Pr(cale_X)}}{\le}
        2^n\exp\bigl(-2(\ln 2)n\bigr)
        = \littleo(1),
      \]
      and the claim holds as required.
    \end{claimproof}

    We conclude showing that $w$ is a Seymour vertex.
    Indeed, since \(w\notin \bad_2\), we have $\deg^+(w)\ge 2\sqrt{\eps}np$.
    Now, Claim~\ref{cl:large-neigh} implies
    that $\bigl|N^2(w)\bigr|\ge (1-2\sqrt{\eps})n$,
    and thus, by Claim~\ref{cl:deg(w)}, we have \(\deg^+(v) < (1-2\sqrt{\eps})n/2\),
    which implies
    \begin{align*}
      \bigl|N_G^2(w)\bigr|
      & \ge   (1-2\sqrt{\eps})n - \deg_G^+(w)
       >     \frac{(1 - 2\sqrt{\eps})n}{2}
        > \deg_G^+(w).\qedhere
    \end{align*}

  \end{proof}

\section{Concluding remarks}\label{sec:concluding-remarks}

In this paper we confirmed Seymour's Second Neighborhood Conjecture (SNC)
for a large family of graphs,
including almost all orientations of (pseudo)random graphs.
We also prove that this conjecture holds a.a.s.\ for
arbitrary orientations of the random graph~$G(n,p)$, where
$p=p(n)$ lies below $1/4$. Interestingly, this range of $p$ encompasses both
sparse and dense random graphs.

The main arguments in our proofs lie in finding a vertex~$w$ of relatively low
outdegree whose outneighborhood contains many vertices of somewhat large outdegree.
Since outneighbors of $w$ cannot have small common outneighborhood,
we conclude that $\bigl|N^2(w)\bigr|$ must be large.

Naturally, it would be interesting to extend further the range of densities for
which arbitrary orientations of $G(n,p)$ satisfy the SNC.

It is seems likely that other classes of graphs,
such as $(n,d,\lambda)$-graphs,
are susceptible to attack using this approach.
Theorem~\ref{t:typical} is also a small step
towards the following weaker
version of Conjecture~\ref{conj:snc}.

\begin{question}
  Do most orientations of an arbitrary graph $G$
  satisfy the SNC?
\end{question}

\phantomsection
\addcontentsline{toc}{section}{Acknowledgments}
\section*{Acknowledgments}

The authors thank Yoshiharu Kohayakawa for useful discussions,
in particular for suggesting we consider bijumbled graphs.

\bigskip

\noindent{\par\small\linespread{1}%\footnotesize
This research has been partially supported by Coordena\c c\~ao de Aperfei\c coamento
de Pessoal de N\'\i vel Superior -- Brasil -- CAPES -- Finance Code 001.
F.~Botler is supported by CNPq {(\small 423395/2018-1)}
and by FAPERJ {(\small 211.305/2019 and 201.334/2022)}.
P.\,Moura is supported by FAPEMIG {(\small APQ-01040-21)}.
T.~Naia is supported by CNPq {(\small 201114/2014-3)} and FAPESP {(\small 2019/04375-5, 2019/13364-7, 2020/16570-4)}.
FAPEMIG, FAPERJ and FAPESP are, respectively, Research Foundations of Minas Gerais,
Rio de~Janeiro and S\~ao Paulo.  CNPq is the National
Council for Scientific and Technological Development of Brazil.

}

% Bibliography -------------------------------------------
\phantomsection
\addcontentsline{toc}{section}{References}
\begin{adjustwidth}{-0em}{-0.1em}
\bibliographystyle{abbrv}
{\footnotesize

  \bibliography{on-gnp.bib}

}
\end{adjustwidth}

\appendix % ---------------------------------------------
\section[Appendix: G(n,p) is p-typical]{Proof that $G(n,p)$ is $p$-typical (Lemma~\ref{l:gnp-typical})}
\label{a:auxiliary}

In this section, we show that $G(n,p)$ satisfies the standard properties
of Definition~\ref{d:p-typical}.
To simplify this exposition, we make use of Lemma~\ref{l:chernoff-2term} below.
Let \defi{$B\sim \calb(N,p)$} denote that $B$ is a binomial
random variable corresponding to the number of successes in
$N$ mutually independent trials,
each with success probability~$p$.

\begin{lemma}\label{l:chernoff-2term}
  For all $N\in\bbn$, all $p\in (0,1)$ and all positive $x$,
  if $B\sim \calb(N,p)$
  then
  \begin{equation*}
    \Pr\bigl[\, |B - Np| > \sqrt{6Np(1-p) x} + 2x  \,\bigr] < 2\exp(-3x).
  \end{equation*}
\end{lemma}

Lemma~\ref{l:chernoff-2term}
follows from the following Chernoff inequality
(see~\cite[Lemma 2.1]{2000:Janson_etal:RG}).

\begin{lemma}\label{l:chernoff}
  Let $X\sim\calb(N,p)$
  and $\sigma^2 = Np(1-p)$.
  For all $t > 0$ we have
  \[
    \Pr\bigl[|X - \EE X | > t\bigr] < 2\exp\left(-\frac{t^2}{2(\sigma^2 + t/3)}\right).
  \]
\end{lemma}

\begin{proof}[Proof of Lemma~\ref{l:chernoff-2term} using Lemma~\ref{l:chernoff}]
  Let $\sigma^2 = Np(1-p)$
  and~\(t = \sqrt{x^2 + 6x\sigma^2} + x\).
  Since~\((t-x)^2 = x^2 + 6x\sigma^2\),
  we have \(t^2 = 2tx + 6x\sigma^2 = 6x(\sigma^2 +t/3)\).
  By Lemma~\ref{l:chernoff},
  \begin{equation}\label{eq:simple-chernoff}
    \Pr\bigl[| B-\EE\, B| > t\bigr]
    < 2\exp\left(-\frac{t^2}{2(\sigma^2 + t/3)}\right) = 2\exp(-3x).
  \end{equation}
  Since
  $t \le \sqrt{6\sigma^2x} + 2x$,
  we have
  \begin{align*}
    \Pr\bigl[| B-\EE\, B| > \sqrt{6\sigma^2x} + 2x\bigr]
    & \le  \Pr\bigl[| B-\EE\, B| > t\bigr]
    \stackrel{\eqref{eq:simple-chernoff}}{<}
       2\exp(-3x). \qedhere
  \end{align*}
\end{proof}

We next show that $G(n,p)$ is $p$-typical.
The properties in Definition~\ref{d:p-typical}
follow by choosing \(x\) in Lemma~\ref{l:chernoff-2term}
so as to make the appropriate a union bound small.

\begin{unnumlemma}[Lemma~\ref{l:gnp-typical}]
  For every $p\colon\mathbb{N} \to (0,1)$,
  a.a.s.\ $G=G(n,p)$ is $p$-typical.
\end{unnumlemma}

\begin{proof}
  We will show that a.a.s.~\ref{i:p-typical:1}--\ref{i:p-typical:common-neigh}
  of Definition~\ref{d:p-typical}
  hold.
  Given a random variable~$Z$ and~$x > 0$, let~$\bbone(Z,x)$ be
  the indicator variable of the `bad' event
  \[
    |Z - \EE Z| > \sqrt{6x\Var(Z)} + 2x,
  \]
  where~$\Var (Z)$ is the variance of~$Z$.
  By Lemma~\ref{l:chernoff-2term}, if $Z\sim\calb(N,p)$
  then
  \begin{equation}\label{e:exp-Z}
    \EE\bigl( \bbone(Z,x)\bigr) = \Pr\bigl[ \bbone(Z,x) = 1\bigr] < 2\exp(-3x).
  \end{equation}

  Firstly, we show that a.a.s.~\ref{i:p-typical:1}
  holds.
  For each $X\subseteq V(G)$, let $Z_X=e(X)$ and let
  \[
    Z^\star=\sum_{X\subseteq V(G)} \bbone\left(Z_X,n\right),
  \]
  taking~$x=n$.
  Note that $Z_X\sim\calb\bigl(\binom{|X|}{2},p\bigr)$ for all $X$.
  By linearity of expectation,
  \begin{align*}
    \EE\, Z^\star
    & = \sum_{X\subseteq V(G)} \EE \bigl(\bbone(Z_X,n)\bigr)
    \stackrel{\eqref{e:exp-Z}}{<}
      \sum_{X\subseteq V(G)} 2\exp\left(-3n\right)
     <
      2^{n+1}\exp\left(-3n\right) = \littleo(1).
  \end{align*}
  Since $Z^\star\ge 0$ (it is the sum of indicator random variables),
  we may use Markov's inequality,
  obtaining~$\Pr[Z^\star\ge 1]\le \EE\,Z^\star = \littleo(1)$.

  A similar calculation, considering in turn
  $\deg(v)$ or~$N(u)\cap N(v)$ instead of~$e(X)$,
  proves that each of the items~\ref{i:p-typical:3} and~\ref{i:p-typical:common-neigh}
  fails to hold with probability~$\littleo(1)$,
  taking~$x$ as~$\ln n$ in both cases, and taking
  union bounds over $n$ or~$\binom{n}{2}$~events respectively.
  Hence $G(n,p)$ satisfies properties~\ref{i:p-typical:1},
  \ref{i:p-typical:3} and~\ref{i:p-typical:common-neigh}
  with probability~$1-\littleo(1)$.

  The strategy to prove~\ref{i:p-typical:sharp-XY}
  is similar to the above, but calculating the number of~events
  in the union bound is slightly more involved.
  If $n'=n''=n$, then (as~above) we consider~$e(X,Y)$ in place of~$e(X)$,
  let $x=n$ and take a union bound over $2^{2n}$ events.
  Otherwise, if $1\le n'\ln n\le n''\le n$,
  then let $\Omega$ be the set of pairs $\{X,Y\}$ with $X,Y,\in V(G)$
  and $|X|, |Y| \leq n'$, and
  note that \(|\Omega| \leq  1+\bigl(\sum_{i=1}^{n'} \binom{n}{i}\bigr)^2\).
  Since $i\le n'<n/2$ for sufficiently large $n$,
  we have \(\binom{n}{i} \leq \binom{n}{n'}\le \left(\frac{en}{n'}\right)^{n'}\)
  and therefore
  \[
    |\Omega|
    \le 1 + \biggl(\:\sum_{i=1}^{n'} \binom{n}{n'}\biggr)^2
    < \left(n'\binom{n}{n'}\right)^2
    < \biggl(\,n' \left(\frac{\ee n}{n'}\right)^{n'}\,\biggr)^2
    < \exp\bigl(2n'(1+\ln n)\bigr)
  \]
  By Lemma~\ref{l:chernoff-2term},
  for each $\{X,Y\}\in \Omega$
  we have~$\Pr\bigl[\bbone(e(X,Y),n'')\bigr] < 2\exp(-3n'')$.
  Applying Markov's inequality
  to~$Z^\star = \sum_{\{X,Y\}\in \Omega}\bbone\bigl(e(X,Y),n''\bigr)$,
  we obtain
  \begin{align*}
    \Pr[Z^\star\ge 1]
    & \le \EE\, Z^\star
      <    \exp\bigl(2n'(1+\ln n)\bigr)\cdot 2\exp(-3n'')
      \le  2\exp(-n''/2)
      =    \littleo(1),
  \end{align*}
  where we  use that $\ln n \le n'\ln n \le n''$.
\end{proof}

\end{document}